\documentclass[a4p,11pt]{article}          % twocolumn
\usepackage[margin=2cm]{geometry}
\usepackage{graphicx}
\usepackage{amsmath,amssymb,bm}
\usepackage[utf8]{inputenc}
\usepackage{microtype}
\usepackage{algorithmic,algorithm}
\usepackage{latexsym}
\usepackage[colorlinks=true,linkcolor=blue,citecolor=blue,urlcolor=blue]{hyperref}
\usepackage[misc,geometry]{ifsym}
\usepackage{subfigure}
\newcommand{\tabref}[1]{Table.~\ref{#1}} 
\newcommand{\figref}[1]{Fig.~\ref{#1}} 
\newcommand{\secref}[1]{Section~\ref{#1}}

\DeclareMathOperator*{\minimize}{minimize}
\begin{document}

\title{Data-driven forced response analysis with min-max representations of nonlinear restoring forces\thanks{
This manuscript is a preprint and has been submitted to a journal for possible publication.
}}
{
%\thanks{Grants or other notes
%about the article that should go on the front page should be
%placed here. General acknowledgments should be placed at the end of the article.}
}
%\subtitle{Do you have a subtitle?\\ If so, write it here}

%\titlerunning{Short form of title}        % if too long for running head

\author{Akira Saito\thanks{Department of Mechanical Engineering, 
Meiji University, Kawasaki, Kanagawa 214-8571, Japan. 
E-mail: asaito@meiji.ac.jp} 
\and 
Hiromu Fujita\thanks{Department of Mechanical Engineering, 
Meiji University, Kawasaki, Kanagawa 214-8571, Japan.}}

\date{}
\maketitle

%\authorrunning{Short form of author list} % if too long for running head

%\institute{A. Saito and H. Fujita\at
%	Department of Mechanical Engineering\\
%	Meiji University\\
%	Kawasaki, Kanagawa 214-8571, Japan\\
%              \email{asaito@meiji.ac.jp}           %  \\
%}

%\date{Received: date / Accepted: date}
% The correct dates will be entered by the editor

%\maketitle

\begin{abstract}
This paper discusses a novel data-driven nonlinearity identification method for mechanical systems with nonlinear restoring forces such as polynomial, piecewise-linear, and general displacement-dependent nonlinearities. The proposed method is built upon the universal approximation theorem that states that a nonlinear function can be approximated by a linear combination of activation functions in artificial neural network framework. The proposed approach utilizes piecewise linear springs with initial gaps to act as the activation functions of the neurons of artificial neural networks. A library of piecewise linear springs with initial gaps are constructed, and the contributions of the springs on the nonlinear restoring force are determined by solving the linear regression problems. The piecewise linear springs are realized by combinations of min and max functions with biases. The proposed method is applied to a Duffing oscillator with cubic stiffness, and a piecewise linear oscillator with a gap and their nonlinearities are successfully determined from their free responses. The obtained models are then used for conducting forced response analysis and the results match well with those of the original system. The method is then applied to experimentally-obtained free response data of a cantilevered plate that is subjected to magnetic restoring force, and successfully finds the piecewise linear representation of the magnetic force. It is also shown that the obtained model is capable of accurately capturing the steady-state response of the system subject to harmonic base excitation . 

% \PACS{PACS code1 \and PACS code2 \and more}
%\subclass{MSC code1 \and MSC code2 \and more}
\end{abstract}
\section{Introduction}
In recent years, data-driven analysis methods have been widely applied in linear and nonlinear mechanical systems for predicting, designing and controlling their dynamics~\cite{BruntonKutz2019}. 

% Koopman or DMD based
Koopman-operator based methods, such as the extended dynamic mode decomposition (DMD) has been applied to nonlinear oscillators~\cite{LiEtAl2017}. Various data-driven DMD-based methods have also been developed to date to predict the dynamics of nonlinear systems. For instance, Galerkin projection method~\cite{AllaKutz2017,SaitoTanaka2023} has been applied to successfully construct reduced order models of nonlinear equations of motion. DMD with control~\cite{ProctorEtAl2016} has also been applied to predict forced response of the nonlinear oscillator~\cite{NamikiSaito2025}. 

% PINNs

% Data-driven approximation of nonlinear forces
A growing number of attempts have been made to represent nonlinear dynamics using artificial neural networks (ANNs). 
For instance, Kuptsov {\it et al.}~\cite{KuptsovEtAl2021} investigated the equivalence of a network model consisting of a single hidden layer of perceptron with sigmoid functions and nonlinear functions in dynamical systems. They examined a couple of nonlinear ODEs including Lorenz system, and showed that the network models exhibit almost equivalent dynamical behavior including bifurcation diagrams and Fourier spectra. 
They also successfully applied ANN to the dynamics of a stiff system~\cite{KuptsovEtAl2023}. 
Moukhliss {\it et al.} proposed an ANN-based surrogate modeling framework for predicting the frequency of nonlinear oscillations of functionally graded materials~\cite{MoukhlissEtAl2026}. 
Luo {\it et al.} proposed a hybrid physics-data-driven method to identify unknown structural nonlinear boundary condition of systems by using both finite element models and multi-layer perceptron (MLP) model~\cite{LuoEtAl2025}. It has been successfully applied to systems with various nonlinearity types including hardening spring, anisotropic spring, and hysteretic nonlinearity. 
These ANN-based nonlinear identification methods are promising because they are versatile and capable of representing a broad class of nonlinear systems. However, the identified model tends to become a black-box and hence it is hard to interpret its physical meanings. 

Raissi {\it et al.} developed a hybrid data-driven and model-based method that achieves system identification and prediction of its dynamics for nonlinear partial differential equations called physics-informed neural networks~\cite{RaissiEtAl2019}. This class of methods enable the solution of nonlinear partial differential equations~(PDEs) by minimizing a cost function that is a combination of the PDE constraint and data errors with respect to the ANN's parameters, which result in the ANNs capable of not only solving the PDEs but also achieving their consistency with the measurement data. This method, however, still result in the black-box nonlinear dynamics model based on the ANNs. 

Brunton {\it et al.} developed a sparse identification of nonlinear dynamical systems~(SINDy)\cite{BruntonEtAl2016}, which utilizes sequential thresholded linear least squares of nonlinear terms by a library of nonlinear functions. 
This gives relatively clear symbolic representation of the nonlinear terms, which tends to provide interpretable nonlinear dynamics models. 

To obtain interpretable nonlinear dynamics model with ANNs, Gonzalez and Lara proposed a nonlinear system identification method for second-order nonlinear ordinary differential equations~(ODEs) based on the concept of characteristic curves where the system dynamics are represented in terms of the nonlinear elements in the ODEs~\cite{GonzalezEtAl2025}. They combined the concept with SINDy or ANN, and the characteristic curve method has been successfully applied to van der Pol equation and mechanical systems with dry-friction. 

Building upon these studies, this paper attempts to propose a simple data-driven system identification and forecasting method for mechanical systems with nonlinear restoring forces with ANN-like system representation combined with a SINDy-like linear regression, which gives us interpretable network model of nonlinear restoring forces as a network of piecewise linear {\it springs}. The proposed approach does not require explicit expressions of the nonlinear forces as in \cite{GonzalezEtAl2025}. 
Fundamentally, the proposed method is based on a piecewise linear approximation of nonlinear functions applied to nonlinear dynamics. 

% Piecewise linear approximation of nonlinear forces
Piecewise linear approximation of nonlinear forces in the dynamics has long been studied by many researchers for various purposes. 
For instance, Garg proposed a piecewise linear approximation of a certain class of nonlinear functions to simplify the numerical integration involved in obtaining the aerodynamic forces and moments from the nonlinear functions~\cite{Garg1981}. 
Cao {\it et al.}~\cite{CaoEtAl2008} proposed a piecewise linear approximation of a smooth discontinuous nonlinear force generated by a linked pair of inclined springs that are pinned to rigid supports. They derived an analytical expression of a trilinear approximation for the nonlinear force using two piecewise linear functions. 
Xie and Shih ~\cite{XieShih2014} also proposed a method to solve for the response of a single degree of freedom (DOF) system with nonlinear stiffness, by utilizing the piecewise linear approximation of nonlinear force and finding the exact closed form solution at each time step. 

More recently, Li {\it et al}~\cite{LiEtAl2024} developed a method for designing a piecewise linear oscillator to mimic nonlinear restoring forces by selecting the breakpoints and slopes of the piecewise linear oscillators. It was shown that the designed piecewise linear oscillator successfully approximates the behavior of a nonlinear system. As it will be shown, the principle behind their method is shared with the proposed approach in this paper. 
Inspired by these studies, the proposed approach in this paper integrates such piecewise linear approximations of nonlinear forces into a single-layer perceptron~(SLP)-like ANN model, where the contributions of the springs are determined from data through SINDy-like linear regression. 

This paper is organized as follows. In \secref{sec:math}, mathematical foundation of the proposed approach is presented. Numerical examples that show the validity of the proposed method are presented in \secref{sec:numerical_example} including a Duffing oscillator and a single DOF system with a piecewise linear stiffness. In \secref{sec:exp}, the proposed method is applied to experimentally-obtained nonlinear vibration data where a plate with a permanent magnet at its tip is subject to nonlinear magnetic force. 
Conclusions of the paper are summarized in \secref{sec:conclusion}. 
%%%%%%%%%%%%%%%%
\section{Mathematical foundation}\label{sec:math}
\subsection{Min-max representation of nonlinear restoring forces}
Let us consider the dynamics of a single mass-spring-damper system of the form, 
\begin{equation}
m\ddot{x}+c\dot{x}+k{x} + f(x) = 0\label{eq1}
\end{equation}
where $m$ is the mass, $c$ is the damping coefficient, $k$ is the spring constant, $x$ is the displacement. Also, let $f(x)$ be a nonlinear function that represents an arbitrary nonlinear force that is a function of $x$. 
\begin{figure}[tb]
\centering
\includegraphics[scale=1]{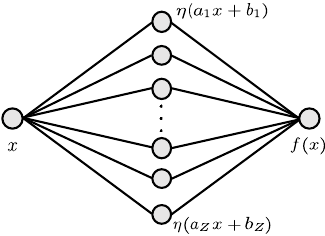}
\caption{Schematic diagram of a single layer perceptron}\label{fig0}
\end{figure}
Based on the universal approximation theorem~\cite{MhaskarMicchelli1992,HornikEtAl1989,Cybenko1989}, which is the foundation of the ANN, assuming that $f(x)$ belongs to a class of tempered distributions, its functional approximation can be found
\begin{equation}
f(x)\approx\sum_{j=1}^{Z}\eta(a_jx+b_j)\label{eq2}
\end{equation}
where $\eta$ are activation functions, $a_j$ are the weights, and $b_j$ are the biases, and $Z$ is the number of neurons in the perceptron.  
The schematic diagram of the SLP is shown in \figref{fig0}. 
The activation function can be simple functions, such as sigmoid function, hyperbolic tangent, and {\it max} function (or so called rectifier linear unit~(ReLU) in ANN). 
The definition of the max function is ${\rm max}(0,x)=x$ for $x>0$ or 0 for $x\leqslant 0$. Also, the {\it min} function can be an activation function where ${\rm min}(0,x)=x$ for $x<0$ or 0 for $x\geqslant 0$. Note that since ${\rm min}(0,x)= -{\rm max}(0,-x)$, they can be related by a linear transformation. 

\begin{figure}[tb]
\centering
\subfigure[Original nonlinear system]{\includegraphics[scale=1]{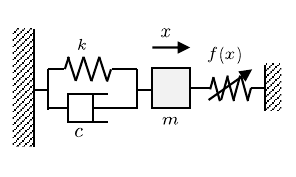}}
\subfigure[Equivalent system representation using a network of piecewise linear springs]{\includegraphics[scale=1]{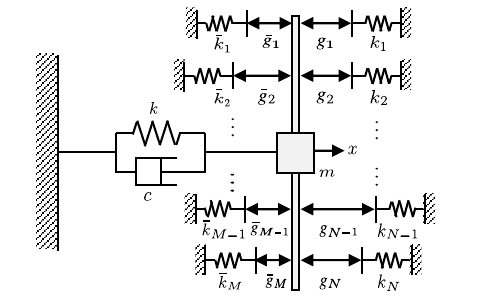}}
\caption{Schematic diagram of multiple piecewise springs acting on the mass}\label{fig1}
\end{figure}

If the min and max functions are chosen as the activation function, $\eta(a_jx+b_j)={\rm max}(0,a_jx+b_j)=a_j{\rm max}(0,x+b_j/a_j)=a_j\eta(x+d_j)$ where $d_j\triangleq b_j/a_j$. 
Similarly, 
$\eta(a_jx+b_j)={\rm min}(0,a_jx+b_j)=a_j{\rm min}(0,x+b_j/a_j)=a_j\eta(x+d_j)$ where $d_j\triangleq b_j/a_j$. 
Defining ${\bf a}\triangleq[a_1,a_2,\dots,a_Z]^{\rm T}$, ${\bf d}\triangleq[d_1,d_2,\dots,d_Z]^{\rm T}$, and $\bm{\eta}(x,{\bf d})\triangleq[\eta(x+d_1),\eta(x+d_2),\dots,\eta(x+d_Z)]^{\rm T}$, then Eq~\eqref{eq2} is rewritten as $f(x)\approx{\bf a}^T\bm{\eta}(x,{\bf d})$. 
With these notations, Eq~\eqref{eq1} becomes
\begin{equation}
m\ddot{x}+c\dot{x}+kx+{\bf a}^{\rm T}\bm{\eta}(x,{\bf d})=0.\label{eq3}
\end{equation}
%%%%%%%%%%%%%
%From the universal approximation theory, following holds:
%\begin{lem}
%Let $f(x)$ be continuous on a compact interval $[x_L,x_R]$. For any $\varepsilon>0$, there exist $M,N\in\mathbb{Z}$, stiffnesses $\bar{k}_j,k_j$, and gaps $\bar{g}_j,\bar{g}_j$, such that
%\begin{equation*}
%{\tiny
%\sup_{x\in[x_L,x_R]}\left|f(x)
%-\sum_{j=1}^M \bar{k}_j\min(0,x-\bar{g}_j)
%-\sum_{j=1}^N k_j\max(0,x-g_j)
%\right|<\varepsilon.}
%\end{equation*}
%\end{lem}
%Therefore, i
If both min and max functions are chosen as the activation functions, Eq.\eqref{eq2} can be interpreted as a linear combination of forces exerted by piecewise linear springs with initial gaps, i.e., 
\begin{equation}
f(x)\approx\sum_{j=1}^{M}\bar{k}_j{\rm min}\left(0,x-\bar{g}_j\right)+
\sum_{j=1}^{N}k_j{\rm max}\left(0,x-g_j\right)
\label{eq4}
\end{equation}
where $\bar{k}_j$ and $k_j$ are the equivalent spring constants of the piecewise linear springs that correspond to $a_j$ of Eq.~\eqref{eq2} for min and max functions respectively, $\bar{g}_j$ and $g_j$ are the equivalent initial gaps that correspond to $\bar{k}_j$ and $k_j$ and $b_j=-\bar{g}_j$ or $b_j=-g_j$. It is noted that both the spring constants and the initial gaps can be negative. 
Theoretically, either max or min function can approximate any nonlinear function. However, to achieve better convergence performance for both increasing and decreasing functions in $x$, both min and max functions are used in the approximation. 
The schematic of this concept is shown in \figref{fig1}. 
Namely, the nonlinear force exerted by the nonlinear stiffness in the original system is replaced by a set of piecewise linear springs with different initial gaps and spring constants. 
This means that {\it any} nonlinear force that is dependent only on $x$ can be modeled as a sum of piecewise linear stiffnesses whose spring constants are $\bar{k}_j$'s and $k_j$'s with the corresponding initial gaps $\bar{g}_j$'s and $g_j$. Namely, the following is equivalent to Eq.\eqref{eq1} for sufficiently large $M$ and $N$, 
\begin{equation}
m\ddot{x}+c\dot{x}+k{x} + 
\sum_{j=1}^{M}\bar{k}_j{\rm min}\left(0,x-\bar{g}_j\right)+
\sum_{j=1}^{N}k_j{\rm max}\left(0,x-g_j\right)
= 0. 
\label{eq5}
\end{equation}
More compactly, 
\begin{equation}
m\ddot{x}+c\dot{x}+k{x} + 
\bar{\bf k}^{\rm T}\bar{\bm{\xi}}(x,\bar{\bf g})+
{\bf k}^{\rm T}{\bm{\xi}}(x,{\bf g})
\label{eq5_2}
= 0, 
\end{equation}
where $\bar{\bm{\xi}}(x,\bar{\bf g})\triangleq[{\rm min}(0,x-\bar{g}_1),\cdots,{\rm min}(0,x-\bar{g}_{M})]^{\rm T}$,
$\bar{\bf k}=[\bar{k}_1,\cdots, \bar{k}_{M}]^{\rm T}$, 
$\bm{\xi}(x,{\bf g})\triangleq[{\rm max}(0,x-{g}_1),\cdots,{\rm max}(0,x-{g}_N)]^{\rm T}$, and ${\bf k}=[k_1,\cdots,k_N]^{\rm T}$. 

The authors have shown that Eq.~\eqref{eq5} holds for a piecewise linear system with an initial gap~\cite{KankiSaito2024}, and the spring constants ${\bf k}$ for the max functions can be successfully obtained from measurement data by using SINDy. 
Namely, a piecewise linear system with an initial gap can be approximated by a sum of piecewise linear systems with initial gaps. 
This paper further explores this concept by examining the equivalence of the system behavior of nonlinear systems with that of the approximated nonlinear system with Eq.~\eqref{eq5_2}. 
\subsection{Determination of equivalent spring constants and gaps}\label{subsec:regression}
In general ANN formation, the weights and biases of the activation functions are determined based on the minimization of a loss function that is typically defined as a squared error between the values of the functions and those of the approximations by optimization solvers such as stochastic gradient methods. 
Since the proposed method uses a linear combination of nonlinear functions of a single layer, simple linear regression of the nonlinear force with a list of nonlinear functions, or a library can be applied as in SINDy. 

Two scenarios are considered. One is the case where the data of the nonlinear function of interest is available and direct fit of the nonlinear function is possible. Such cases are discussed in \ref{direct_method} and referred to as direct method. The other scenario is where the data of the nonlinear function of interest is not available but measurements data of $x$ is available. This case is discussed in \ref{indirect_method} and referred to as indirect method. 
\subsubsection{Direct method}\label{direct_method}
Assuming that the oscillator moves within a certain range, e.g., $x_{L}\leqslant x\leqslant x_{R}$, and the nonlinear restoring force $f(x)$ is measured at finite and discrete measurement points $\tilde{\bf x}=[x_1,\cdots,x_n]^{\rm T}$, 
i.e., ${\bf f}(\tilde{\bf x})=[f(x_1),\cdots,f(x_n)]^{\rm T}$ 
where $x_{L}=x_1$ and $x_{R}=x_n$ and $x_1<\cdots<x_n$. Then, the equivalent spring constants of ${\rm min}$ and ${\rm max}$ functions in Eq.~\eqref{eq4} can be obtained simply by solving the following linear least squares problem. Namely, we seek $\tilde{\bf k}$ and ${\bf k}$ such that 
\begin{equation}
\minimize_{\bm{\kappa}}: \lVert\bm{\mathcal{L}}(\tilde{\bf x},\bar{\bf g},{\bf g})\bm{\kappa}-{\bf f}(\tilde{\bf x})\rVert^2, 
\label{eq_minim}
\end{equation}
%which requires, 
%\begin{equation}
%\bm{\mathcal{L}}(\tilde{\bf x},\bar{\bf g},{\bf g})\bm{\kappa}={\bf f}(\tilde{\bf x}), 
%\end{equation}
where $\bm{\kappa}=[\bar{\bf k}^{\rm T}, {\bf k}^{\rm T}]^{\rm T}$ contains a list of unknown equivalent spring constants to be determined. 
$\bar{\bf g}$ and ${\bf g}$ contain lists of $M$ and $N$ possible gap values that need to be set a priori. 
$\bm{\mathcal{L}}(\tilde{\bf x},\bar{\bf g}, {\bf g})$ contains a library of min and max functions evaluated at $x_1,\cdots,x_n$, i.e., $\bm{\mathcal{L}}(\tilde{\bf x},\bar{\bf g},{\bf g})=[\bar{\bf L}(\tilde{\bf x},\bar{\bf g}),{\bf L}(\tilde{\bf x},{\bf g})]$, where 
\begin{equation}
\bar{\bf L}(\tilde{\bf x},\bar{\bf g})=
\begin{bmatrix}
{\rm min}(0,x_1-\bar{g}_1), \cdots,{\rm min}(0,x_1-\bar{g}_{M})\\
\vdots\\
{\rm min}(0,x_n-\bar{g}_1), \cdots,{\rm min}(0,x_n-\bar{g}_{M})
\end{bmatrix}
\end{equation}
and 
\begin{equation}
{\bf L}(\tilde{\bf x},{\bf g})=
\begin{bmatrix}
{\rm max}(0,x_1-{g}_1), \cdots,{\rm max}(0,x_1-{g}_{{N}})\\
\vdots\\
{\rm max}(0,x_n-{g}_1), \cdots,{\rm max}(0,x_n-{g}_{M})
\end{bmatrix}.
\end{equation}
The minimizer of Eq.~\eqref{eq_minim} can be obtained as 
\begin{equation}
\bm{\kappa}=\left[\bm{\mathcal{L}}(\tilde{\bf x},\bar{\bf g},{\bf g})\right]^{\dagger}{\bf f}(\tilde{\bf x}), \label{eq:kstar}
\end{equation}
where $^{\dagger}$ denotes the Moore-Penrose pseudo-inverse of a rectangular matrix. Obviously, this process of obtaining the contributions from the min and max functions does not require sophisticated gradient-based optimization methods such as stochastic gradient descent or Adam~\cite{KingmaBa2014} that are used in the construction of modern ANNs. This comes at the expense of forming a list of possible gaps ($\bar{g}_i$'s and $g_i$'s) that correspond to the biases of the neurons in ANNs, which are typically obtained by the optimization solvers. Instead, the contributions from the piecewise linear springs that are in the list ($\bar{\bf L}$ and ${\bf L}$) are obtained as the spring constants ($\bm{\kappa}$), which correspond to the weights of the activation functions of the neurons, by the regression Eq.~\eqref{eq:kstar}. 
As long as the moving range of the oscillator can be estimated a priori, the list of gaps can be created relatively easily from the moving range.  Unnecessary gap values are then naturally eliminated in the process of linear regression. 
\subsubsection{Indirect method}\label{indirect_method}
This case also assumes that the oscillator moves within a certain range $x_L\leqslant x\leqslant x_R$. 
This time, it is assumed that the direct measurement of nonlinear force is not feasible and hence its data is unavailable. Instead, displacement, velocity, and acceleration of the nonlinear system are measured at discrete time instants $t=[t_1,t_2,\cdots t_n]$, and they are used to indirectly identify the nonlinear restoring force. For this purpose, regression of time series of the available data is conducted as follows. 
Rearranging the terms in Eq.~\eqref{eq5_2}, we obtain, 
\begin{equation}
\ddot{x}+2\zeta\omega_{\rm n}\dot{x}+\omega_{\rm n}^2 x =- \bar{\bf k}'^{\rm T}\bar{\bm{\xi}}(x,\bar{\bf g})-{\bf k}'^{\rm T}\bm{\xi}(x,{\bf g}), 
\label{eq:aaa}
\end{equation}
where $\zeta=c/2\sqrt{mk}$, $\omega_{\rm n}=\sqrt{k/m}$, $\bar{\bf k}'=\bar{\bf k}/m$, and ${\bf k}'={\bf k}/m$. Note that $\omega_{\rm n}$ is the natural frequency and $\zeta$ is the damping ratio with the absence of the nonlinear terms on the right hand side of Eq.~\eqref{eq:aaa}. 
Let 
\begin{align*}
\tilde{\bf a}&=[\ddot{x}(t_1),\cdots,\ddot{x}(t_n)]^{\rm T}=[a_1,\cdots,a_n]^{\rm T}, \\
\tilde{\bf v}&=[\dot{x}(t_1),\cdots,\dot{x}(t_n)]^{\rm T}=[v_1,\cdots,v_n]^{\rm T}, \\
\tilde{\bf x}&=[x(t_1),\cdots,x(t_n)]^{\rm T}=[x_1,\cdots,x_n]^{\rm T}. 
\end{align*}
We then seek $\bar{\bf k}$ and ${\bf k}$ such that
\begin{equation}
\minimize_{\bm{\kappa}'}:\lVert
\mathcal{\bm{L}}(\tilde{\bf x},\bar{\bf g},{\bf g})\bm{\kappa}'-\tilde{\bf a}_\ell
\rVert^2,\label{eq_indirect}
\end{equation}
where $\bm{\kappa}'=[\bar{{\bf k}}'^{\rm T},{{\bf k}}'^{\rm T}]^{\rm T}$, and 
$\tilde{\bf a}_\ell$ corresponds to the linear portion of the equation of motion evaluated with the responses of the nonlinear system, or 
$\tilde{\bf a}_\ell\triangleq\tilde{\bf a}+2\zeta\omega_{\rm n}\tilde{\bf v}+\omega_{\rm n}^2\tilde{\bf x}$. 
Therefore, the minimizer of Eq.~\eqref{eq_indirect} is obtained as, 
\begin{equation}
\bm{\kappa}'=\left[\bm{\mathcal{L}}
(\tilde{\bf x},\bar{\bf g},{\bf g})
\right]^{\dagger}\bar{\bf a}_{\ell}.
\end{equation}
With the obtained $\bm{\kappa}'$, the right hand side of Eq.~\eqref{eq:aaa} can be evaluated for any $x$. 
\section{Numerical examples}\label{sec:numerical_example}
\subsection{Duffing oscillator}
\begin{figure}[tb]
\centering
\includegraphics[width=9cm]{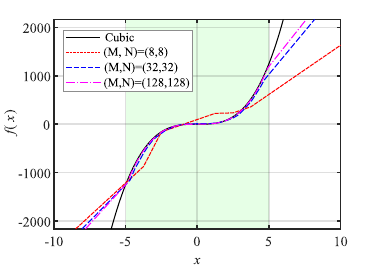}
\caption{Convergence of the approximated function with respect to the number of functions. $p_2=10$. The regression was conducted for data in the filled region ($-5\leqslant x \leqslant 5$).}
\label{fig2}
\end{figure}

\begin{figure}[tb]
\centering
\subfigure[$(M,N)=(8,8)$]{\includegraphics[width=8.5cm]{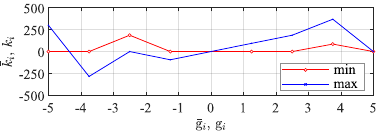}}
\subfigure[$(M,N)=(32,32)$]{\includegraphics[width=8.5cm]{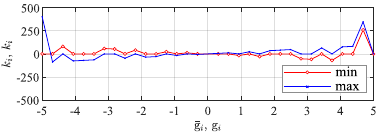}}
\subfigure[$(M,N)=(128,128)$]{\includegraphics[width=8.5cm]{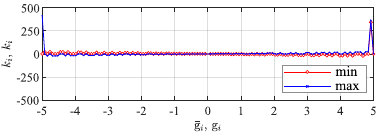}}
\caption{Spring constants versus the initial gap values}\label{fig2_11}
\end{figure}

The proposed method is applied to a single DOF system with a cubic restoring force, or a Duffing oscillator of the form:
\begin{equation}
m\ddot{x}+c\dot{x}+p_1{x}+p_2x^3=0. \label{sec2:eq1}
\end{equation}
where $p_1$ and $p_2$ are arbitrarily chosen stiffness parameters. 
First, the equivalent form of the nonlinear term $p_2x^3$ represented by Eq.~\eqref{eq4} needs to be found. 
Namely, the spring constants $k_j$ need to be computed, or they need to be trained using known data of $f(x)$. 
In this case, the spring constants can be obtained by simple regression using the direct method as described in \ref{direct_method}. 
\subsubsection{Piecewise linear approximation of the cubic term}
Figure \ref{fig2} shows the results of the approximations using the min-max functions for $(M,N)=(8,8),(32,32)$, and $(128,128)$ where $p_2=10$.  The spring constants $\bar{\bf k}$ and ${\bf k}$ were obtained for the values in $-5\leqslant x\leqslant5$. As can be seen, the approximation is poor for $(M,N)=(8,8)$. 
However, as $M$ and $N$ increase, the approximation becomes more accurate. Indeed, when $(M,N)=(128,128)$, as can be seen, the original function is quite well represented by the piecewise linear functions especially within $-5\leqslant x\leqslant5$. The approximation quality is degraded outside the range of $-5\leqslant x\leqslant5$, but it is expected because the spring constants $\bar{\bf k}$ and ${\bf k}$ were trained by the dataset in $-5\leqslant x\leqslant5$ and the capability of the approximation to extrapolate the function outside the range is not guaranteed.  

In addition, \figref{fig2_11} shows the obtained equivalent spring constants versus the gap values for different values of $M$ and $N$. As seen in \figref{fig2_11}(a), the springs constants are large for large initial gaps near the lower and upper bound of the range. This makes sense because the cubic nonlinearity produces the largest values at the bounds and hence the springs with large spring constants with large initial gaps need to be added to represent such a nonlinear force. 

\subsubsection{Free response}
\begin{figure}[tb]
\centering
\includegraphics[width=8.5cm]{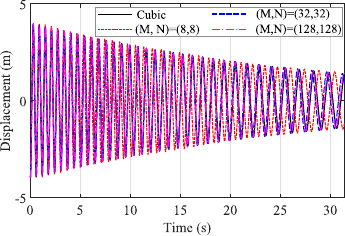}
\caption{Time history of the oscillator with the original cubic nonlinearity and with the min-max representations}
\label{sec1:fig1}
\end{figure}
\begin{figure}[tb]
\centering
\subfigure[With original cubic stiffness]{\includegraphics[width=8cm]{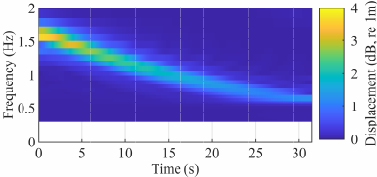}}
\subfigure[With $(M,N)=(8,8)$]{\includegraphics[width=8cm]{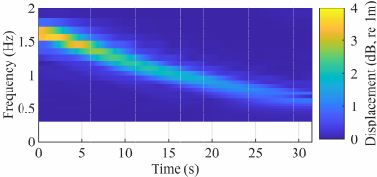}}
\subfigure[With $(M,N)=(128,128)$]{\includegraphics[width=8cm]{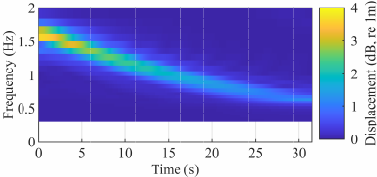}}
\caption{Comparison of the time-frequency relationships using scalogram of the response}
\label{sec1:fig2}
\end{figure}
\begin{table}[tb]
\centering
\caption{RMSE values for the approximations}\label{sec1:tab1}
\begin{tabular}{cccc}\\ \hline \hline
$(M,N)$ & $(8,8)$ & $(32,32)$ & $(128,128)$\\ \hline
RMSE (m) & $2.38$ & $0.497$ & $0.215$ \\ \hline
\end{tabular}
\end{table}
The original equation of motion Eq.~\eqref{sec2:eq1}, and the ones with the nonlinear term replaced with the min-max representations with the trained $\bar{\bf k}$ and ${\bf k}$ have been solved by time integration for the initial condition of $x(0)=-4.0$, and $\dot{x}(0)=0$. The time histories are shown in \figref{sec1:fig1}. 
As can be seen, the predicted time history of $x(t)$ with all approximation cases match well with the original data until approximately $t=5$s. However, prediction accuracy gradually decreases, and the delay from the original data becomes visible, especially for $(M,N)=(8,8)$. On the other hand, the one with $(M,N)=(32,32)$ does not show significant delay until around $t=25$s. The approximation with $(M,N)=(128,128)$ shows almost negligible delay before $t=30$s. 
\tabref{sec1:tab1} shows the root mean square error (RMSE) values between the approximations and the original response. The RMSE value decreases as $(M,N)$ increases. This also supports the convergence of the proposed approach with respect to $M$ and $N$. 

To see the frequency contents in the response more clearly, wavelet transform has been applied to the original response, the ones with $(M,N)=(8,8)$, and $(M,N)=(128,128)$. The obtained scalograms are shown in \figref{sec1:fig2}. Overall, the trend of the response of the original system can be well captured by the approximations. However, for the case with $(M,N)=(8,8)$, the dominant frequency components near $t=30$s become different from those of the original system. On the other hand, we can see that the frequency contents can be well captured by the approximation with $(M,N)=(128,128)$. 

\subsubsection{Forced response}
\begin{figure}[tb]
\centering
\includegraphics[width=8.5cm]{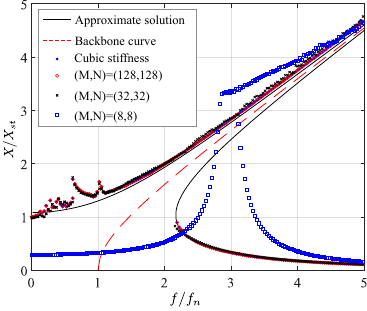}
\caption{Comparison of the forced response results with single cubic stiffness and multiple PWL springs, along with classical approximate solutions and backbone curve}\label{fig:duffing_harmonic}
\end{figure}
To examine the applicability of the proposed approach to forced response problems, a forcing term has been added to Eq.~\eqref{sec2:eq1}, i.e., 
\begin{equation}
m\ddot{x}+c\dot{x}+p_1x+p_2x^3=F\sin(2\pi f t), 
\label{eq_duffing_forced}
\end{equation}
where $F$ is the amplitude of forcing, $f$ is the excitation frequency. Also, the same forcing term is added to the equivalent equation of motion with the cubic term replaced by its equivalent min-max representation, i.e.,
\begin{equation}
m\ddot{x}+c\dot{x}+p_1x+\bar{\bf k}^{\rm T}\bar{\bm{\xi}}(x,\bar{\bf g})+{\bf k}^{\rm T}{\bm{\xi}}(x,{\bf g})=F\sin(2\pi f t). 
\label{eq_duffing_forced_equiv}
\end{equation}
The steady-state solutions of the equations Eqs.~\eqref{eq_duffing_forced} and \eqref{eq_duffing_forced_equiv} have been sought for various frequencies over $0\leqslant{f}\leqslant{5}{\mathrm Hz}$ by using time integration for 200 cycles of a period of excitation. The results are shown in \figref{fig:duffing_harmonic} for the steady-state responses of Eqs.~\eqref{eq_duffing_forced} and \eqref{eq_duffing_forced_equiv} for $(M,N)=(8,8)$, $(32,32)$, and $(128,128)$. When solving the equations for the steady-state, the excitation frequency $f$ was swept upward or downward to obtain solutions of multiple stable branches. 
Backbone curve for the undamped system ($c=0$) and the approximate solutions by a classical perturbation method~\cite{Meirovitch2001} have also been plotted in \figref{fig:duffing_harmonic}. 
Note that $f$ is normalized with respect to the linear natural frequency $f_n=\sqrt{k/m}/(2\pi)$, and the amplitude of oscillation is normalized with respect to the static equilibrium $X_{st}$ that satisfies $p_1X_{st}+p_2X_{st}^3-F=0$. Note that it was numerically obtained as $X_{st}=0.3930$. 
As can be seen, the response curve of the original system with the cubic stiffness shows a typical stiffening effect where the response curve is bent toward high frequency and shows the primary resonance. Also, there are secondary resonances below $f/f_n=1$ including 3:2 superharmonic resonance. We can now see that these characteristics are well captured by the approximate system with $(M,N)=(128,128)$, including the primary and the secondary resonances. On the other hand, for the case with $(M,N)=(32,32)$, the values tend to be slightly off from those of the original system or the approximate system with $(M,N)=(128,128)$. Furthermore, if the numbers of min-max functions are further reduced down to $(M,N)=(8,8)$, the system fails to predict the accurate steady-state response for almost the entire frequency range considered. 
Interestingly, it shows a strongly nonlinear yet typical response curve of piecewise linear systems~(see, e.g., \cite{TienEtAl2023}) at around $X/X_{st}=3.2$, i.e., $X\approx1.25$. This approximately corresponds to the smallest gap value of the active min and max functions, i.e., the system acts almost as a piecewise linear system with the smallest gap value, which does not well represent the Duffing oscillator. 
%%%%%%%%%%%%%%%%%%%%%%%%%%%%%%%%%%%%%%%%%%%%%%%%%%%%%%5
\subsection{Piecewise linear system}
\begin{figure}[tb]
\centering
\includegraphics[width=6cm]{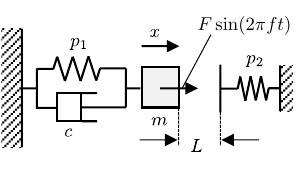}
\caption{Schematic diagram of the piecewise linear system}\label{fig:bilinear_schematic}
\end{figure}
Next, a single DOF system with a piecewise linear stiffness subjected to harmonic forcing is considered. 
Such a piecewise linear system has been studied extensively due to its simple yet strong nonlinearity~\cite{ShawHolmes1983}. 
The schematic diagram of the system is shown in \figref{fig:bilinear_schematic}. The piecewise linear spring is fixed at the right end and the initial gap between the spring and the mass is $L$. 
The equation of motion of the system is written as, 
\begin{equation}
m\ddot{x}+c\dot{x}+p_1x+p_2\mathrm{max}(0,x-L)=F\sin(2\pi f_e t),
\label{eq1:pwl}
\end{equation}
where $p_2$ is the spring constant of the piecewise linear stiffness. The piecewise linear spring with gap $L$ is approximated using the min and max functions with different gap values (without $L$) using the procedure described in \ref{subsec:regression}. Note that the equivalence of a single max function with a gap and multiple max functions with different gap lengths was considered in Ref.~\cite{KankiSaito2024} both numerically and experimentally, and free response was well captured by the max functions. 
\begin{figure}[tb]
\centering
\includegraphics[scale=1]{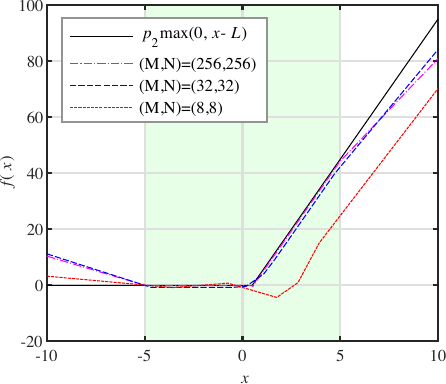}
\caption{Convergence of the approximated function with respect to the number of functions. $p_2 = 10$. The regression was conducted for data in the filled region ($-5\leqslant x \leqslant 5$).}\label{fig:binear_f_vs_x}
\label{fig:bilinear_approximation}
\end{figure}
\begin{figure}[tb]
\centering
\subfigure[$(M,N)=(8,8)$]{\includegraphics[width=8.5cm]{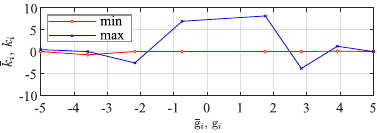}}
\subfigure[$(M,N)=(32,32)$]{\includegraphics[width=8.5cm]{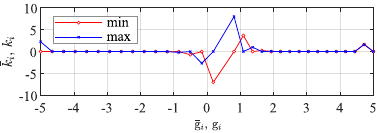}}
\subfigure[$(M,N)=(256,256)$]{\includegraphics[width=8.5cm]{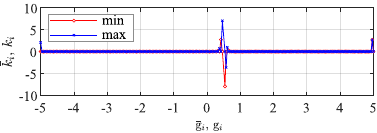}}
\caption{Spring constants versus the initial gap values for the piecewise linear stiffness ($L=0.5$)}\label{fig2_1}\label{fig:bilinear_k_vs_g}
\end{figure}
With the assumption that the data of nonlinear force is available, direct method described in \ref{direct_method} can be applied. 
For this numerical experiment, $L=0.5$ is assumed. 
The linear regression was conducted using data taken within $-5\leqslant x\leqslant 5$ for $(M,N)=(8,8)$, $(32,32)$, and $(256,256)$. 
The approximated functions as well as the original piecewise linear stiffness are shown in \figref{fig:bilinear_approximation}. 
As can be seen, the approximation is poor for $(M,N)=(8,8)$ even within $-5\leqslant x\leqslant 5$. On the other hand, for $(M,N)=(32,32)$, the min and max functions well represent the original function. With $(M,N)=(256,256)$, the difference is negligible for $-5\leqslant x\leqslant 5$. 

The obtained spring constants for the corresponding gap values are plotted in \figref{fig:bilinear_k_vs_g}. As can be seen, even though the original function consists only of the max function, both min and max functions contribute to the representation of a single max function. As the number of min and max functions increases, we can see that the contributions of the min and max functions whose gap values are close to the gap value of the original max function ($L$) become dominant. This appears as a spike in the plot as shown in \figref{fig:bilinear_k_vs_g}(c). This characteristics can be exploited to identify unknown $L$, as it has already been discussed in Ref.~\cite{KankiSaito2024}. 
\subsection{Free response}
\begin{figure}[tb]
\centering
\includegraphics[width=7cm]{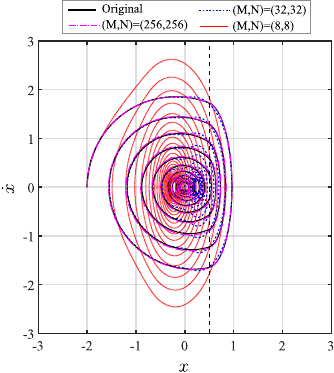}
\caption{Comparison of the free response results with a single piecewise linear stiffness and multiple PWL stiffnesses.}
\label{fig:bilinear_free}
\end{figure}
\begin{table}[tb]
\centering
\caption{RMSE values for the approximations}\label{sec1:tab2}
\begin{tabular}{cccc}\\ \hline \hline
$(M,N)$ & $(8,8)$ & $(32,32)$ & $(256,256)$\\ \hline
RMSE (m) & $0.554$ & $0.229$ & $0.028$ \\ \hline
\end{tabular}
\end{table}

Free response analyses have been conducted for the original piecewise linear system and the system with the multiple min and max functions for $x(0)=-2$ and $\dot{x}(0)=0$. Namely, $F=0$ in Eq.~\eqref{eq1:pwl}. Other parameter values are set to $m=1$, $c=0.1$, and $k=1$. 
The results are shown in \figref{fig:bilinear_free} as phase portrait format to clearly visualize the changes in the dynamics as the system stiffness changes. The vertical dashed line designates $L=0.5$ where the mass collides with the piecewise linear stiffness. The phase portrait with $(M,N)=(8,8)$ poorly predicts the correct trajectory. The one with $(M,N)=(32,32)$ predicts that of the original piecewise linear system well, but there are noticeable differences especially near the gap. On the other hand, $(M,N)=(256,256)$ predicts quite well the one with the original system even near the gap region. Table~\ref{sec1:tab2} shows the RMSE values for the approximations, where the reduction in RMSE can be observed for increasing $M$ and $N$. 
This indicates the high accuracy of the proposed system representation of the nonlinear function even when the nonlinearity is piecewise linear. 
\subsection{Forced response}
\begin{figure}[tb]
\centering
\includegraphics[width=8.5cm]{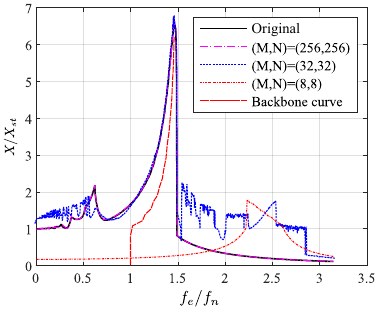}
\caption{Comparison of the forced response results with a single piecewise linear stiffness and multiple PWL springs}\label{fig:bilinear_forced}
\end{figure}
Forced response analyses have been conducted for $0\leqslant f_e\leqslant 0.5$, or $0\leqslant f_e/f_n\leqslant\pi$ where $f_n=1/2\pi$. At each frequency, the equation of motions were solved numerically for $x(0)=0$ and $\dot{x}(0)=0$ until they reach the steady-state solutions. The amplitude of the response was computed and normalized with respect to the static displacement.  It is plotted for the normalized excitation frequency and shown as \figref{fig:bilinear_forced} for the original piecewise linear system, the approximate systems with $(M,N)=(8,8)$, $(32,32)$, and $(256,256)$. As can be seen, there is a primary resonance at $f_e/f_n=1.47$, whereas the largest secondary resonance can be seen at $f_e/f_n=0.622$.  The proposed model with $(M,N)=(8,8)$ fails to predict both primary and secondary resonances. The proposed model with $(M,N)=(32,32)$ succeeded in predicting the primary resonance, however, it failed to capture the secondary resonance and responses at higher frequency range than the primary resonance. The proposed model with $(M,N)=(256,256)$ captures both the primary and the secondary resonances. 
In addition to the forced response results, backbone curve is shown in \figref{fig:bilinear_forced}. The backbone curve was extracted from the instantaneous frequency and amplitude of the free response of the system with the identified multiple PWL stiffnesses using the initial condition of $x(0)=2.0$ and $\dot{x}(0)=0$. The instantaneous amplitude of the response was defined as the mean value of the upper and lower envelope of the response. The instantaneous frequency of the response was computed from the time instants at which the displacement crossed the instantaneous center of oscillation. The instantaneous center of oscillation was defined as the midpoint between the instantaneous upper and lower envelopes. As seen, the backbone curve of the obtained system captures the amplitude dependence of the resonant frequency quite well. 

\begin{figure}[tb]
\centering
\subfigure[Primary resonance: $f_{e}/f_n=1.470$]{\includegraphics[scale=1]{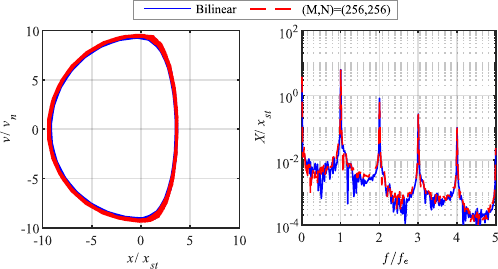}}
\subfigure[Secondary resonance: $f_{e}/f_n=0.622$]{\includegraphics[scale=1]{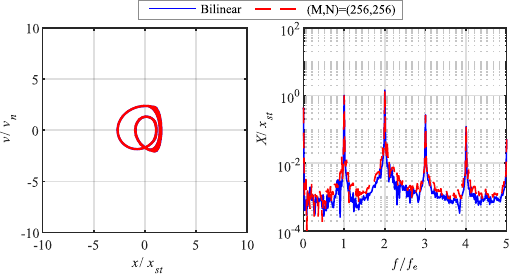}}
\caption{Phase portraits and Fourier spectra at primary and secondary resonances}
\label{fig:harmonic_pwl_spectra}
\end{figure}

To further explore the validity of the proposed approach, trajectory of the solutions as well as their Fourier spectra were examined. Figure \ref{fig:harmonic_pwl_spectra} shows the phase portraits and the Fourier spectra computed by Fast Fourier Transform (FFT). 
The horizontal axis of the phase portrait is normalized with respect to the static displacement $x_{st}$ and the vertical axis is normalized with respect to $v_n=\omega_n x_{st}$ where $\omega_n=\sqrt{k/m}=1$. The horizontal axes of the FFT spectra are normalized with respect to the frequency of the external forcing $f_e$, whereas the vertical axes are normalized with respect to $x_{st}$. 
At the primary resonance shown in \figref{fig:harmonic_pwl_spectra}(a), the response is dominated by the harmonic component of the excitation frequency. Other frequency components that are integer multiples of the excitation frequency are much smaller than the one with $f/f_e=1$. On the other hand, at the secondary resonance, the frequency component with twice the excitation frequency is as large as that with the excitation frequency, as can be seen in \figref{fig:harmonic_pwl_spectra}(b). 
This indicates the existence of the superharmonic motion. 
Indeed, the loop in the phase portrait has a knot, indicating the superharmonic motion at twice the excitation frequency. 
This arises because the excitation frequency approximately equals to twice the natural frequency of the system, resulting in 2:1 superharmonic resonance. 

\section{Experimental validation}\label{sec:exp}
In this section, the proposed approach is applied to a cantilevered plate with magnetic-force induced nonlinearity from experimentally obtained dataset and its applicability to real measurement data is discussed. 
First, the nonlinear restoring force acting on the system is identified by the proposed approach from free response data. Second, forced response of the system under base excitation is predicted by adding the excitation term in the obtained equation of motion and its validity is discussed by comparing it with the measurement. 
\subsection{Experimental setup}
\begin{figure}[tb]
\centering
\includegraphics[scale=1]{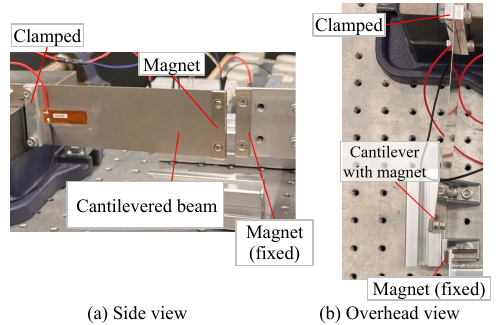}
\caption{Test equipment for free vibration test. The plate is in equilibrium state. }\label{fig:exp:1}
\end{figure}
The system studied in this section is shown in \figref{fig:exp:1}. A cantilevered plate with a permanent magnet at its tip is clamped and fixed to a surface plate. The plate is made from Titanium Alloy and its dimension is $175\mbox{mm}\times50\mbox{mm}\times0.5\mbox{mm}$. 
There is another permanent magnet that is fixed to an aluminum block rigidly attached to the surface plate, and is placed near the tip of the cantilevered plate. 
The distance between the magnet attached to the cantilevered plate and the other magnet was set to 12 mm. 
The directions of the polarity of the magnets are both in the direction perpendicular to the surface of the plate, but in the opposite directions with each other. 
The magnet attached to the cantilevered plate experiences repulsive magnetic force if its position is not aligned with the fixed magnet, as shown in \figref{fig:exp:1}(b). Since the plate also experiences the restoring force due to the elasticity of the plate, there is a stable equilibrium position where the restoring elastic force balances the repulsive force. The same applies to the situation when the plate is bent toward the opposite side. Therefore, this oscillator has two stable equilibria and hence becomes a bistable oscillator. This class of oscillator has been used as a type of energy harvesting devices~\cite{DengEtAl2019} and nonlinear energy sinks~\cite{LiEtAl2026}. 

The expected single DOF approximation of the equation of motion of the system is given as 
\begin{equation}
\ddot{x}+2\zeta\omega_{\rm n}\dot{x}+\omega_{\rm n}^2{x}=f_t(x)\label{exp:eq:0}
\end{equation}
where $x$ is the displacement of the plate at the tip, $f_t(x)$ is the nonlinear restoring force that stem from the magnetic repulsive force generated between the fixed magnet and the magnet at the tip of the plate. Note that the equivalent stiffness coefficient $k$ mainly represents the linear restoring force that comes from the elasticity of the plate. 

We consider the identification of such nonlinear restoring force exerted by the magnets, and aim to develop mathematical model that can well represent the dynamics of the system from measured data. In this case,  the magnetic nonlinear force is not directly measured. Instead, displacement is measured by laser displacement sensor (LDS), and the indirect method introduced in \secref{indirect_method} is applied. 

\subsection{Potential constraints on the basis functions}

For the identification of the magnetic force, the following constraints on the $\mathrm{min}$ and $\mathrm{max}$ functions are considered. 
%Since the magnetic force is known to be dependent on the position of the magnet, it is natural to define its potential. 
The magnetic force acting between the magnets is dependent on the angle, alignment, and the relative distance of the magnets. 
In the present setup, the amplitude of vibration is sufficiently small so that the changes in the angle and the alignment are negligible. Therefore, it is reasonable to assume that the magnetic force is dependent solely on the relative distance of the magnets. Since one of the magnets is fixed to the base, the relative distance between the magnets is dependent only on the displacement of the magnet that is attached to the tip of the cantilever plate. As a consequence, the magnetic force can be regarded as a conservative force, and associated potential function can be defined. 
Namely, we hypothesize that the potential can be approximated as a linear combination of a linear function, quadratic function, and piecewise quadratic functions of the form: 
\begin{equation}
V_t(x,{\bf g})=V_c(x) + V_\ell(x) + V_p(x,{\bf g})\label{exp:eq:1_3}
\end{equation}
where $V_c(x)= q_1 x$, $V_\ell(x)=\frac{1}{2}q_2 x^2$, 
$V_p(x,{\bf g})=\bm{\kappa}^{\rm T}\bm{\phi}(x,{\bf g})$, 
$\bm{\kappa}=[\kappa_1,\dots,\kappa_N]^{\rm T}$, 
$\bm{\phi}(x,{\bf g})=[\phi(x,g_1),\dots,\phi(x,g_N)]$, and
\begin{equation}
\phi(x,g_i)\triangleq 
\frac{1}{2}x^2 + \frac{1}{2}{\rm min}(0,x+{g}_i)^2+\frac{1}{2}{\rm max}(0,x-g_i)^2. 
\label{exp:eq:1}
\end{equation}
$q_1$, $q_2$, and $\kappa_i's$ are the contributions from these terms determined by the regression. 
The nonlinear force is then obtained as:
\begin{align}
f_t(x,{\bf g})&=-\frac{{\rm d}V_t}{{\rm d}x}(x,{\bf g}) \nonumber \\
&=f_c(x)+f_\ell(x)+f_p(x,{\bf g})\label{exp:eq:1_2}
\end{align}
where $f_c(x)=-q_1$, $f_\ell(x)=-q_2x$, $f_p =-\bm{\kappa}^{\rm T}{\bf f}(x,{\bf g})$, ${\bf f}(x,{\bf g})=[f(x,g_1),\dots,f(x,g_N)]^{\rm T}$, and 
\begin{align}
f_i(x,g_i)&=- \frac{{\rm d}}{{\rm d}x}\phi(x,g_i)\nonumber \\
&=- \left\{x+{\rm min}(0,x+{g}_i)+{\rm max}(0,x-g_i)\right\}.\label{exp:eq:2}
\end{align}
The expression on the right hand side is defined as, 
\begin{equation}
\psi(x,g_i)\triangleq \frac{{\rm d}}{{\rm d}x}\phi(x,g_i)=
x+{\rm min}(0,x+{g}_i)+{\rm max}(0,x-g_i)\label{exp:eq:2_1}
\end{equation}
and $\psi$ is herein referred to as basis functions and used for the formation of the library. 
With these assumptions, the physics behind the data can be enforced as constraints on the regression, i.e., potential function associated with the identified force can be defined, otherwise the min-max functions do not respect the physics and hence may produce models that do not follow physical laws. 
\subsection{Free response}

With the system shown in \figref{fig:exp:1}, free response experiment has been conducted. 
The procedure of the experiment is stated as follows. 
The tip of the plate was displaced along transverse direction. The tip was then released with zero initial velocity. 
The displacement of the plate tip where the magnet is attached was measured by an LDS (IL-100, KEYENCE, Japan). The sampling rate was 1kHz. 
To obtain the velocity, derivative of the obtained discrete time history of the displacement was taken numerically. The derivative of the velocity was then taken to obtain the acceleration. Note that a low-pass filter with cut-off frequency of 100Hz has been applied to the measured displacement to eliminate excessive noise in the derivatives of the displacement.  
\begin{figure}[tb]
\centering
\subfigure[Free response of the plate tip]{\includegraphics[width=8.5cm]{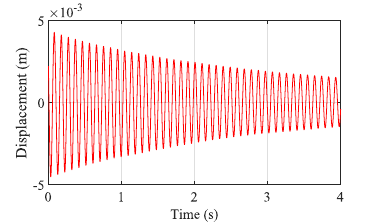}}
\subfigure[Backbone curve obtained from the free response]{\includegraphics[width=8.5cm]{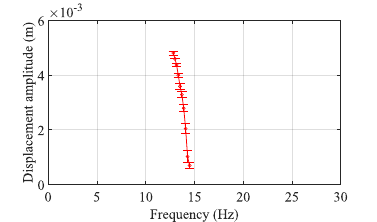}}
\caption{Time response of the displacement of the plate tip and backbone curve obtained from the free response}\label{fig:exp:2}
\end{figure}

Result of the free response measurement is shown in \figref{fig:exp:2}(a). The plate tip vibrates with approximately 13Hz. The amplitude of the displacement gradually decrease due to damping. As a result, the frequency of oscillation slightly deviates as the amplitude of oscillation decreases. To further examine this, the backbone curve was extracted from the data of \figref{fig:exp:2}(a) by the instantaneous frequency and amplitude, and is shown in \figref{fig:exp:2}(b). As can be seen, clear softening effect can be observed. Namely, the backbone curve is bent toward left, which means that the resonant frequency decreases as the amplitude of oscillation increases, or the resonant frequency increases as the amplitude of oscillation decreases. This makes sense because if the tip displacement amplitude increases, the magnetic repulsive force weakens because the distance between the pole of the magnet at the tip and that of the other magnet fixed at the base increases. 
\begin{figure}[tb]
\centering
\subfigure[Potential]{\includegraphics[width=8.5cm]{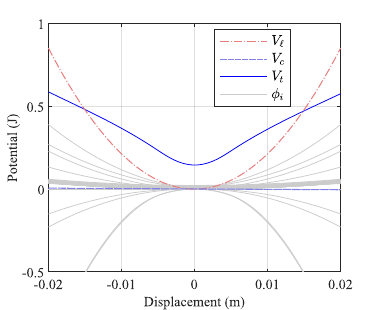}}
\subfigure[Magnetic force]{\includegraphics[width=8.5cm]{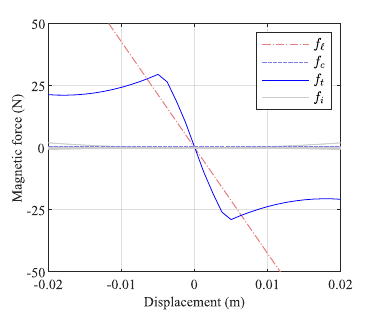}}
\caption{Identified potential and its associated magnetic force}\label{fig:exp:3}
\end{figure}

With the free response shown in \figref{fig:exp:2}, the proposed method has been applied. 

First, the damping ratio $\zeta$ and the natural frequency $\omega_{\rm n}$ corresponding to $m$, $c$ and $k$ in Eq.~\eqref{exp:eq:0} were identified from a free response of the same plate without the magnetic force, i.e., fixed magnet was removed from the fixture. 
The identified values were $\omega_{\rm n}= 65.2\mathrm{rad/s}$, which correspond to 10.4Hz, and $\zeta = 0.0054$. 

Second, the indirect method has been applied on the free response of the system with the magnetic force. 
The library used for the regression of the nonlinear force consists of a constant term 
%, a linear term, 
and the piecewise linear forces defined in Eq.~\eqref{exp:eq:1_2}, i.e., 
\begin{equation}
\bm{\mathcal{L}}(\tilde{\bf x},{\bf g})=
\begin{bmatrix}
1, &\psi(x_1,g_1),&\dots,&\psi(x_1,g_M)\\
& & \vdots& \\
1, &\psi(x_n,g_1),&\dots,&\psi(x_n,g_M)\\
\end{bmatrix}. 
\end{equation}
where $g_1,\dots,g_M$ were equally sampled $M$ points within $-0.02\leqslant x\leqslant 0.02$. 
The identified potential function and the associated magnetic force are shown in \figref{fig:exp:3}. 
In \figref{fig:exp:3}(a), the obtained total potential energy $V_t(x)$ is shown along with the constituent potential energies with respect to displacement $x$, where $V_c(x)$ is the potential corresponding to the constant force, $V_\ell(x)$ is the potential corresponding to linear force, $\phi_i(x)$'s are the piecewise quadratic potential defined as Eq.~\eqref{exp:eq:1_3}. 
As seen in \figref{fig:exp:3}, $V_t(x)$ is a convex function with two inflection points that are almost symmetric with respect to the vertical axis. 
In \figref{fig:exp:3}(b), the obtained total magnetic force $f_t(x)$ is shown with the constant term $f_c(x)$, linear term $f_\ell(x)$, and piecewise linear terms $f_i(x)$. 
As can be seen in the figure, the magnetic force is almost symmetric with respect to the origin. The slope of $f_t$ is slightly steeper than the linear portion of the force $f_\ell$ due to the contributions from quadratic and piecewise quadratic terms. 
At the inflection points of the potential $V_t(x)$, $f_t(x)$ takes its extrema. 

\begin{figure}[tb]
\centering
\subfigure[$M=20$]{\includegraphics[scale=1]{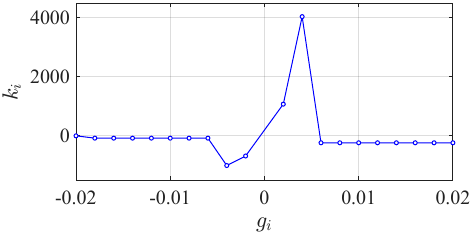}}
\subfigure[$M=32$]{\includegraphics[scale=1]{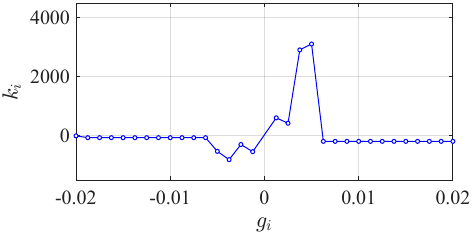}}
\subfigure[$M=40$]{\includegraphics[scale=1]{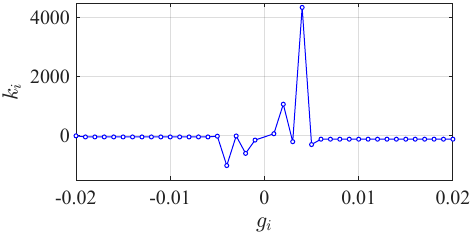}}
\caption{Spring constants versus the initial gap values for the min-max functions}\label{fig:exp:3_1}
\end{figure}
Figure \ref{fig:exp:3_1} shows the identified spring constants, or the contributions of the basis functions $\psi(x,g)$ with respect to $g_i$ for $M=20$, 32, and 40. As $M$ increases, the positive contribution near $g_i=0.005$ tends to become the largest whereas the negative contribution near $g_i=-0.005$ tends to become the smallest. This generally implies that the effects of the piecewise linear springs tend to become larger as the distance from the origin increases. The springs with positive gaps tend to have positive spring constants whereas those with negative gaps tend to have negative spring constants.  The combined effect of these springs represent the observed softening nature of the nonlinearity. 
%This trend is similar to the one observed for the equivalent spring constants found for the Duffing oscillator in \secref{sec:numerical_example}. 
\begin{figure}[tb]
\centering
\includegraphics[width=8.5cm]{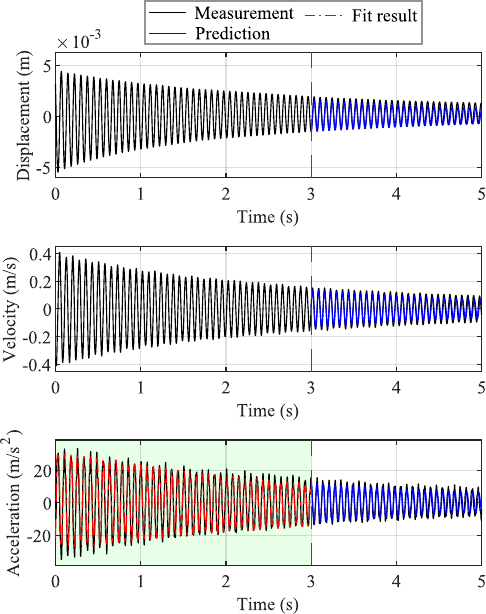}
\caption{Results of free response ($M=32$). The regression was conducted for data in the filled region ($0\leqslant t \leqslant 3$).}\label{fig:exp:4}
\end{figure}

The time histories of the measured displacement, which is the same data as \figref{fig:exp:2}(a), velocity, and acceleration are shown in \figref{fig:exp:4} for $0\mathrm{s}\leqslant t \leqslant 5\mathrm{s}$. 
The data for $0\mathrm{s}\leqslant t \leqslant 3\mathrm{s}$ was used for obtaining magnetic restoring force shown in \figref{fig:exp:3}(b) by the regression, and the resulting nonlinear governing equations were solved by numerical time integration to forecast the time histories of displacement, velocity, and acceleration for $3\mathrm{s}\leqslant t \leqslant 5\mathrm{s}$. As seen, the forecasted displacement, velocity, and acceleration agree very well with the measurements. 
\begin{figure}[tb]
\centering
\includegraphics[scale=1]{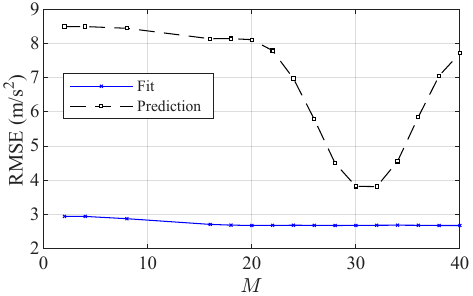}
\caption{Convergence study with respect to the number of $\phi(x,g_i)$ (Fit duration=50\%)}\label{fig:exp:4_1}
\end{figure}

\begin{figure}[tb]
\centering
\includegraphics[scale=1]{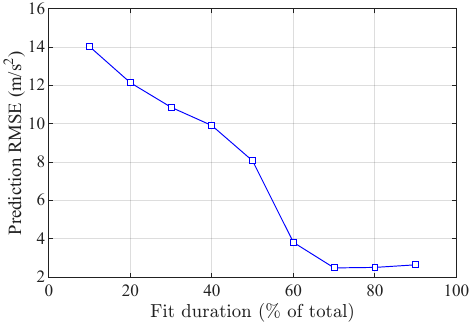}
\caption{Convergence study with respect to the fit duration to total duration ($M=32$)}\label{fig:exp:4_2}
\end{figure}

To further examine the characteristics of the obtained model, convergence study with respect to the number of $\phi(x,g_i)$, or $M$, has been conducted. 
The RMSE value between the measured acceleration and the fit result was computed for $0\leqslant t\leqslant 3$. Also, the RMSE between the measured and the predicted values for $3\leqslant t \leqslant 5$ was computed. This means that the data $0\leqslant t\leqslant 3$ was used for training the model, and that for $3\leqslant t\leqslant 5$ was used for validation. 
The results are shown in \figref{fig:exp:4_1}. The RMSE value for the fit decreases monotonically with respect to $M$. 
On the other hand, the RMSE for the prediction decreases as $M$ increases up to around $30$. For $M\geqslant 30$, the RMSE for the prediction increases. 
Considering that the RMSE for the fit decreases for $M\geqslant 30$, this clearly shows the consequence of overfitting.  
This shows that convergence study needs to be conducted with respect to $M$, or cross-validation needs to be conducted to find the optimal $M$. 

Furthermore, convergence study with respect to the duration of the fit to the total time has also been conducted. Fit duration was varied from 0.5s to 4.5s while the total time was fixed to 5.0s. Namely, the fit duration was varied from 9.98\% to 90.0\%. $M$ was fixed to 32. The prediction RMSE decreases as fit duration increases up to 70\%. It then slightly increases for fit duration larger than 70\%. This trend implies that even though the prediction accuracy increases as fit duration increases for moderately small fit duration, if the fit duration exceeds 70\%, it suffers from overfitting. This also shows that the convergence study with respect to the prediction accuracy should also be conducted to obtain the optimal fit duration. 
\subsection{Forced response}
\begin{figure}[tb]
\centering
\includegraphics[scale=1]{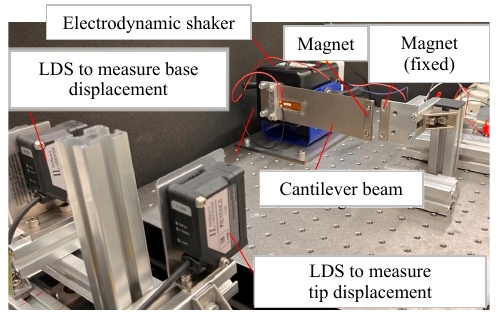}
\caption{Test equipment for forced vibration test}\label{fig:exp:5}
\end{figure}
Next, forced response of the oscillator is examined. Namely, the prediction capability of the model obtained in the previous section is examined under harmonic excitation. The obtained results are compared with measurements. 

\begin{figure}[tb]
\centering
\subfigure[Experiment (phase shift: 178.20deg]{\includegraphics[scale=1]{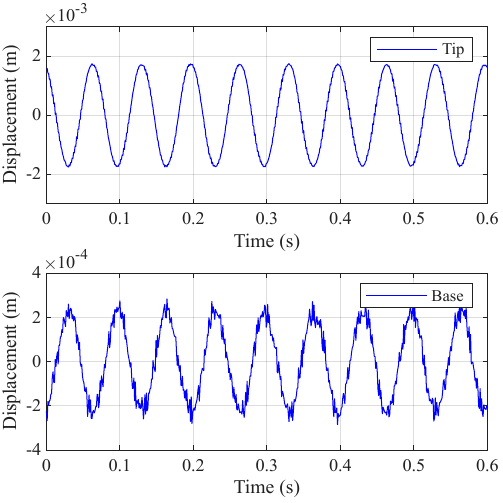}}
\subfigure[Simulation (phase shift: 175.86deg) ]{\includegraphics[scale=1]{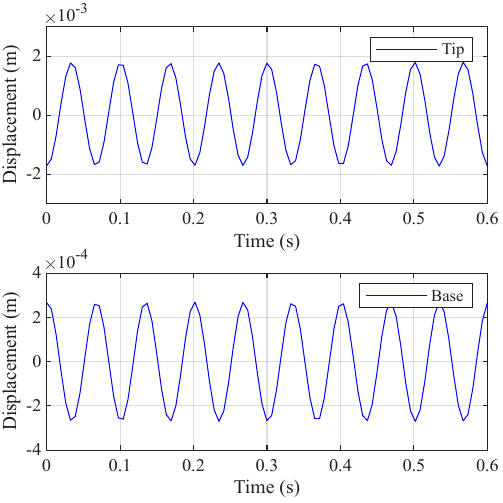}}
\caption{Comparison of the measured forced response time histories and those computed with the proposed approach at 15Hz. }\label{fig:exp:6}
\end{figure}

The forced response test was conducted by using an electrodynamic shaker. The test setup is shown in \figref{fig:exp:5}. The oscillator was fixed at the moving head of the electrodynamic shaker (K2004E01, The Modal shop, USA). A sinusoidal signal was supplied from a function generator (WF1947, NF Corporation, Japan) to the electrodynamic shaker to generate the base excitation to the oscillator. The displacement was measured at the moving head of the electrodynamic shaker where the plate is attached (herein referred to as base) and at the tip of the plate (herein referred to as tip) using two of the LDSs. 
With this setup, forced response tests have been conducted. The sinusoidal excitation with a fixed frequency has been applied to the plate, and its vibration response was measured at the tip and the base until it reaches the steady-state. The measurement was repeated for different frequencies from 5Hz to 30Hz with 1Hz increment. 
To examine the amplitude-dependence of the system, the shaker was operated with two vibration amplitude levels, i.e., 0.27mm and 0.83mm, which approximately correspond to the input DC voltages to the shaker of 0.1V and 0.3V, respectively. 
As a representative result, \figref{fig:exp:6}(a) shows the results  of forced response at the steady-state for the sinusoidal excitation with its frequency of 15Hz. The displacements measured at the tip and the base are shown. 
The response was extracted from the measurements for an arbitrary time window of 0.6s in the steady-state. 
As can be seen, the sinusoidal motion of the base induces almost sinusoidal motion of the tip with phase shift. The phase shift between the base and the tip motions were computed to be 178.20deg. 

Forced response calculation has been conducted with the model identified from the free response. The equation of motion Eq.~\eqref{exp:eq:0}, which is unforced, was converted to that with the base excitation, i.e., 
\begin{align}
\ddot{x}+2\zeta\omega_{\rm n}\dot{y}+\omega_{\rm n}^2y&=f_t(y,{\bf g})\\
\mathrm{with}\quad y&\triangleq x-X_b\sin (2\pi f_e t)
\end{align}
where $y$ is the relative displacement of the oscillator measured from the base, $X_b$ is the amplitude of the base excitation, $f_e$ is the excitation frequency of the base. For the evaluation of the nonlinear force $f_t(y,{\bf g})$, the parameters $q_1$, $q_2$, and $\bm{\kappa}$ in Eq~\eqref{exp:eq:1_2} obtained from the free response analysis were used. 

Figure \ref{fig:exp:6}(b) shows the steady-state response predicted by the model generated based on the free-response. The response was extracted from the results of numerical integration for an arbitrary time window of 0.6s in the steady-state. 
As can be seen, the tip displacement prediction shown in \figref{fig:exp:6}(b) agree well with the measurement shown in \figref{fig:exp:6}(a), in terms of the amplitude of oscillation and the phase shift from the base displacement. The phase shift was computed to be 175.86deg, which is 1.3\% error from the measurement. 

\begin{figure}[tb]
\centering
\includegraphics[scale=1]{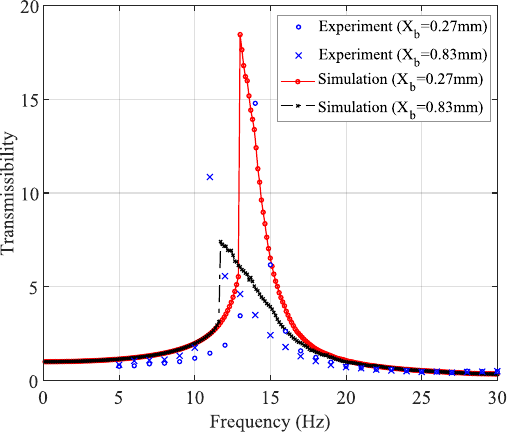}
\caption{Comparison of the measured forced response results and those computed with the proposed approach}\label{fig:exp:7}
\end{figure}

Figure \ref{fig:exp:7} shows the results of forced response near the resonance in terms of transmissibility and frequency for both experiments and simulations. The transmissibility was defined as the ratio of the steady-state amplitude of oscillation of the tip to that of the base ($X_b$). 
The plots clearly show that the system has softening nonlinearity, i.e., the resonant peak bends toward lower frequency. 
The graph shows that the resonant frequency is lower for $X_b=0.83$mm than that for $X_b=0.27$mm. On the other hand, resonant peak amplitude is higher for $X_b=0.83$mm than that for $X_b=0.27$mm. 
%This makes sense because in general the resonant frequency of a system with softening nonlinearity decreases as the amplitude of oscillation increases. {\color{red}{Check Nayfeh \& Mook}}
These trends are well captured by the model obtained by the proposed approach. 
\section{Conclusions}\label{sec:conclusion}
In this paper, a novel data-driven nonlinear identification method for mechanical systems under nonlinear restoring forces has been proposed. The method is based upon the network of piecewise linear springs represented by min and max functions. The contributions from the springs on the nonlinear restoring forces are obtained by solving a regression problem of the nonlinear forces with a library of min-max functions. 
The method has been applied to a Duffing oscillator and a piecewise linear oscillator. It was shown that the resulting models can capture the free and forced response under harmonic excitations. The method was also applied to identify the magnetic restoring forces exerted on a cantilevered plate. To exploit the conservative nature of the magnetic forces, potential constraints are applied to the min-max representation of the nonlinear force such that they represent conservative forces. The proposed method successfully identified the nonlinear restoring forces exerted by the permanent magnet on the cantilevered plate and the associated potential function and produced a model that can accurately predict free and steady-state forced response of the system under harmonic excitation. 

%\begin{acknowledgements}
%
%\end{acknowledgements}
%
%
%% Authors must disclose all relationships or interests that 
%% could have direct or potential influence or impart bias on 
%% the work: 
%%
%% \section*{Conflict of interest}
%%
%% The authors declare that they have no conflict of interest.
%
%
% BibTeX users please use one of
%\bibliographystyle{spbasic}      % basic style, author-year citations
\bibliographystyle{unsrt}
%\bibliographystyle{spmpsci}      % mathematics and physical sciences
%\bibliographystyle{spphys}       % APS-like style for physics
%\bibliography{./references.bib}   % name your BibTeX data base

\section*{Statements and Declarations}

\noindent{\bf Competing interests} The authors have no relevant financial or non-financial interest to disclose. \\ \\
\noindent{\bf Author contributions} 
A. Saito: Conceptualization, Methodology, Formal analysis, Investigation (numerical), Writing -- original draft, Writing -- review \& editing. H. Fujita: Resources, Investigation (experimental), Data curation, Validation, Writing -- review \& editing. All authors approved the final manuscript. 
\\ \\
\noindent{\bf Data availability} The dataset generated during and/or analyzed during the current study are available from the corresponding author upon reasonable request. \\ \\
%% Non-BibTeX users please use
%\begin{thebibliography}{}
%%
%% and use \bibitem to create references. Consult the Instructions
%% for authors for reference list style.
%%
%\bibitem{RefJ}
%% Format for Journal Reference
%Author, Article title, Journal, Volume, page numbers (year)
%% Format for books
%\bibitem{RefB}
%Author, Book title, page numbers. Publisher, place (year)
%% etc
%\end{thebibliography}

\end{document}